\documentclass[12pt]{amsart}

\usepackage{amsfonts}
\usepackage{dur}

\mathchardef\Gamma="100
\mathchardef\Delta="101
\mathchardef\Pi="105

\def\C{{\mathbb C}}
\def\fm{\varphi}

\def\gkn{G^k(\C^n)}

\def\m{\psi}
\newcommand{\Mat}{\operatorname{Mat}}
\def\mnk{\Mat(n, k; \C)}
\def\mnkt{\Mat(n, k; \C)^\times}
\def\mbar#1{\kern 0.1em\overline{\kern -0.1em #1 \kern -0.1em} 
  \kern 0.1em}
\def\pr{\qopname \relax o{pr}}
\newcommand{\rank}{\operatorname{rank}}
\def\scf{w}
\def\sm{\xi}

\newcommand{\tr}{\operatorname{tr}}

\newtheorem{prop}{Proposition}[section]
\newtheorem{theor}{Theorem}[section]

\numberwithin{equation}{section}

\begin{document}

\title{Frenet frames and Toda systems}

\author{A. V. Razumov}

\address{Institute for High Energy Physics, 142284 Protvino, Moscow Region,
Russia. \vskip 0pt \tt E-mail: razumov@mx.ihep.su}

\begin{abstract}
It is shown that the integrability conditions of the equations satisfied by
the local Frenet frame associated with a holomorphic curve in a complex
Grassmann manifold coincide with a special class of nonabelian Toda
equations. A local moving frame of a holomorphic immersion of a Riemann
surface into a complex Grassmann manifold is constructed and the
corresponding connection coefficients are calculated. 
\end{abstract}

\maketitle

\section{Introduction}

In the present paper we consider local differential geometry of a
holomorphic curve $\m$ in the complex Grassmann manifold $\gkn$. It is well
known that the simplest way to study such a curve is to consider a
holomorphic local lift $\sm$ of $\m$ to the space $\mnk$ of complex $n
\times k$ matrices, defined on some open subset $U$ of $M$. The natural
object arising here is the so called Frenet frame being a mapping which
associates a basis of $\C^n$ with a point of $U$. We give a procedure of
the construction of the Frenet frame and derive the equations which are
satisfied by the corresponding $\C^n$-valued functions on $M$. Then we show
that the obtained equations are closely connected to nonabelian Toda
systems based on the Lie group ${\rm GL}(n, \C)$. The equations describing
these systems are exactly integrable, so we get a possibility to construct
holomorphic embeddings with prescribed properties. Note that our definition
of the Frenet frame is different from the definition usually used for the
problem under consideration, see, for example, \cite{Li97}. Namely, we do
not perform the complete orthogonalisation of the basis, and it is our
definition that leads to an exactly integrable system. The paper is
concluded with the explanation of how the Frenet frame can be used to
construct a local moving frame associated with the mapping $\m$ and with
the calculation of the corresponding connection coefficients.

There is a lot of papers discussing geometry of harmonic mappings of a
Riemann surface into a projective space, see, for example, \cite{EWo83,
BWo92} and references therein. Since any holomorphic mapping from one
K\"ahler manifold to another one is harmonic, most results of those paper
are valid for holomorphic mappings. There are also some papers devoted to
consideration of harmonic mappings from Riemann surfaces to Grassmann
manifolds \cite{EWo83, CWo85, BWo86, BSa87, CWo87, Wol88, Woo88, Val88,
Uhl89}. In these papers the problems of classification and construction of
harmonic mappings are mainly investigated.  On the other hand, local
differential geometry properties of holomorphic mappings from Riemann
surfaces to Grassmann manifolds have not yet been discussed in a systematic
way. In this respect we would like to mention the paper \cite{Li97} which
deals with differential geometry of holomorphic mappings from the
two-sphere to a Grassmann manifold.

Note that the connection of abelian Toda systems based on the Lie group
${\rm SL}(n, \C)$ with the equations satisfied by the Frenet frame
associated with an immersion of a Riemann surface into a projective space
appeared to be useful for the study of $W$-geometries \cite{GMa92, GMa93,
Ger93}. Moreover, this connection allows to give a simple proof of the
infinitesimal Pl\"ucker relations. The abelian Toda systems based on other
semisimple Lie groups can be associated with equations satisfied by the
Frenet frames of holomorphic curves in the corresponding flag manifolds.
This leads to the so-called generalised infinitesimal Pl\"ucker relations
\cite{RSa94, GSa96}. The results of the present paper show that there are
some relations between the geometric characteristics of the holomorphic
curves associated with a holomorphic curve in a Grassmann manifold, but
here we have no direct analogues of the Pl\"ucker relations. This problem
requires an additional consideration. It is also very interesting to
investigate the corresponding $W$-geometries.

\section{Frenet frames}

\subsection{Holomorphic curves in Grassmann manifolds} \label{s:2}

We start with the description of the realisation of complex Grassmann
manifolds which is convenient for our purposes. Let $\mnk$ denote the
complex manifold of $n \times k$ complex matrices. Assume that $n > k$ and
denote by $\mnk^\times$ the open submanifold of $\mnk$ formed by $n \times
k$ matrices of rank $k$. The manifold $\mnk^\times$ could be naturally
identified with the set of all bases of $k$-dimensional subspaces of
$\C^n$. Two points $p, p' \in \mnk^\times$ correspond to bases of the same
subspace of $\C^n$ if and only if there is an element $g$ of the complex
general linear group ${\rm GL}(k, \C)$ such that $p' = pg$. Thus we have
the free right action of the complex Lie group ${\rm GL}(k, \C)$ on
$\mnk^\times$ connecting different bases of the $k$-dimensional subspaces
of $\C^n$. The orbit space here is the Grassmann manifold $\gkn$. The
canonical projection $\pi: \mnk^\times \to \gkn$ is holomorphic and
$\mnk^\times \stackrel{\pi}{\to} \gkn$ is a holomorphic principal ${\rm
GL}(k, \C)$-bundle.

Consider a holomorphic curve $\m : M \to \gkn$, where $M$ is some connected
Riemann surface. Any such a curve defines a holomorphic subbundle
$\underline \m$ of rank $k$ of the trivial holomorphic fibre bundle
$\underline\C^n = M \times \C^n$. Here the fiber over the point $p \in M$
is $\m(p)$. The mapping $\m$ is said to be {\it linearly full} if the
bundle $\underline \m$ is not contained in any proper trivial subbundle of
$\underline \C^n$.

To describe geometry of the corresponding holomorphic curve in $\gkn$ it is
convenient to use local holomorphic lifts of the mapping $\m$ to $\Mat(n ,
k; \C)^\times$. We say that a holomorphic mapping $\sm: U \to \mnk^\times$,
where $U$ is an open subset of $M$, is a local holomorphic lift of $\m$ if
it satisfies the relation
\[
\pi \circ \sm = \m|_U.
\]
Since our consideration is local, we will use local holomorphic lifts
defined in coordinate neighbourhoods of $M$. The corresponding local
coordinate will be denoted by $z$. Moreover, to make our notations
consistent with those used in description of Toda systems, we denote $z^+ =
z$, $z^- = \bar z$ and
\[
\partial_- = \partial/\partial z, \qquad \partial_+ = \partial/\partial
\bar z.
\]

Let $\sm: U \to \mnk^\times$ be a lift of a holomorphic mapping $\m: M \to
\gkn$. Denote by $f_1$, $\ldots$, $f_k$ the holomorphic mappings from $U$
to $\C^n$ determined by the columns of the matrix-valued mapping $\sm$.
For any $p \in U$ the linear subspace $W_p$ of $\C^n$ spanned by the
vectors $f_1(p)$, $\ldots$, $f_k(p)$ and by the vectors given by the first
and higher order derivatives of the functions $f_1$, $\ldots$, $f_k$ over
$z$ at the point $p$, is independent of the choice of the lift $\sm$ and
local coordinate $z$.  

\begin{prop} \label{t:2}
For any points $p$ and $p'$ of $M$ one has
\[
W_p = W_{p'}.
\]
\end{prop}

\begin{proof}
Let $p$ be an arbitrary point of $U$ and $c = z(p)$. Expand the functions
\begin{equation}
f_1 \circ z^{-1}, \; \ldots, \; f_k \circ z^{-1} \label{**}
\end{equation}
in the power series at $c$, which converges absolutely in some disc $D(c,
r) \subset z(U)$. Since the first and higher order derivatives of the
functions  $f_1 \circ z^{-1}$, $\ldots$, $f_k \circ z^{-1}$ over $z$ can be
also expanded at $c$ in the power series which absolutely converge in $D(c,
r_c)$, it is clear that for any $p' \in z^{-1}(D(c, r))$ one has $W_{p'}
\subset W_p$. From the other hand, if $p' \in z^{-1}(D(c, r/2))$ then the
functions (\ref{**}) and their first and higher derivatives over $z$ can be
expanded at the point $c' = z(p')$ in a power series which converge at
least in the disc $D(c', r/2)$. Since $p \in z^{-1}(D(c', r/2))$, one
concludes that $W_p \subset W_{p'}$. Therefore, for all points $p'$ such
that $|z(p) - z(p')| < r/2$ one has $W_{p'} = W_p$. 

Considering various local lifts of $\m$ we conclude that for any point $p
\in M$ there is an open neighbourhood $U_p$ such that $W_{p'} = W_p$ for
all $p' \in U_p$. Now the assertion of the Proposition follows from the
connectedness of $M$.
\end{proof}

Thus we have some constant subspace $W$ of $\C^n$ which characterises the
curve $\m$.  The mapping $\m$ can be considered as a mapping from $M$ to
$G^k(W)$, in this case it can be easily shown that $\m$ is linearly full.
Therefore, without any loss of generality we can restrict ourselves to the
consideration of linearly full mappings. 

\subsection{Some properties of vector-valued holomorphic functions}
\label{s:1}

In this Section $V$ is an $n$-dimensional complex vector space and $U$ is
an open connected subset of the complex plane $\C$. The results similar to
ones discussed below in this Section for the case of polynomial mappings
are proved in \cite{Li95, Li97}.

A mapping from $U$ to $V$ is called a {\it $V$-valued function} on $U$. Let
$f$ be a $V$-valued function on $U$. For any $w \in V^*$ one defines the
complex function $w(f)$ by
\[
w(f)(c) = w(f(c)).
\]
If $g$ is a $V^*$-valued function on $U$ one defines the complex function
$g(f)$ by
\[
g(f)(c) = g(c)(f(c)).
\]
The function $f$ is said to be {\it holomorphic\/}, if for any $w \in V^*$
the function $w(f)$ is holomorphic. 

Let $\{e_i\}$ be a basis of $V$, $\{e^i\}$ be the dual basis of $V^*$, and
$f$ be a $V$-valued function on $U$. It is easy to verify that for any $c
\in U$ one has
\[
f(c) = \sum_{i=1}^n e_i \, b^i(c),
\]
where $b^i = e^i(f)$. Here it is customary to write
\begin{equation}
f = \sum_{i=1}^n e_i \, b^i. \label{1}
\end{equation}
It is clear that the function $f$ is holomorphic if and only if the
functions $b^i$ are holomorphic.

If $f$ is a holomorphic $V$-valued function on $U$ and $g$ is a holomorphic
$V^*$-valued function on $U$, then $g(f)$ is a holomorphic function on $U$.

A point $c \in U$ is called a {\it zero} of a holomorphic $V$-valued
function $f  \not\equiv 0$ on $U$, if $f(c) = 0$. A zero of a holomorphic
$V$-valued function $f \not\equiv 0$ is a common zero of the functions
$b^i$ entering a representation of type (\ref{1}). Here, in general, the
order of the zero is different for different values of the index $i$,
moreover, it depends on the choice of the basis $\{e_i\}$. Nevertheless,
the minimal value of order does not depend on the choice of the basis
$\{e_i\}$ and it characterises the $V$-valued function $f$ itself. The
corresponding positive integer is called the {\it order} of zero of $f$. 

A holomorphic $V$-valued function $f \not\equiv 0$ on $U$ may have only a
countable number of isolated zeros, and the set of zeros has no limit
points in $U$.

If $c \in U$ is a zero of order $\mu$ of a holomorphic $V$-valued function
$f \not\equiv 0$ on $U$, one can represent $f$ as
\[
f = g (z-c)^\mu.
\]
Here $g$ is a holomorphic $V$-valued function on $U$ having the same zeros
as $f$, except the zero at $a$.

\begin{prop} \label{p:1}
Let $f \not\equiv 0$ be a holomorphic $V$-valued function on $U$. There
exists a holomorphic $V$-valued function $g$ on $U$ having no zeros in $U$
and such that there is valid the representation
\[
f = g \, d,
\]
where $d$ is a holomorphic function on $U$.
\end{prop} 

\begin{proof}
In accordance with the generalisation of the Weierstrass theorem
\cite{Mar77}, there exists a holomorphic function $d$ having the same zeros
and with the same orders as $f$. It is easy to show that the function $g =
f/d$ satisfies the requirement of the Proposition.
\end{proof}

The {\it rank} of a set $\{f_1, \ldots, f_k\}$ of holomorphic $V$-valued
functions on $U$ is defined as
\[
\rank \{f_1, \ldots, f_k\} = \max_{c \in U} \rank \{f_1(c), \ldots,
f_k(c)\}.
\]
The rank of a set $\{f_1, \ldots, f_k\}$ of holomorphic $V$-valued
functions on $U$ equals $k$ if and only if $f_1 \wedge \cdots \wedge f_k
\not\equiv 0$. We say that a set $\{f_1, \ldots, f_k\}$ of holomorphic
$V$-valued functions on $U$ is of {\it constant rank} if
\[
\rank \{f_1(c), \ldots, f_k(c)\} = \rank \{f_1, \ldots, f_k\}
\]
for any $c \in U$. 

\begin{prop} \label{p:2}
Let the set $\{f_1, \ldots, f_k\}$ of holomorphic $V$-valued functions on
$U$ be of rank $k$. Then there exists a set $\{g_1, \ldots, g_k\}$ of
$V$-valued holomorphic functions on $U$ of constant rank $k$ such that 
\[
f_\alpha = \sum_{\beta=1}^k g_\beta \, d^{\beta}_\alpha, \qquad \alpha = 1,
\ldots, k,
\]
where $d^\beta_\alpha$ are holomorphic functions on $U$.
\end{prop} 

\begin{proof}
We prove the Proposition by induction over $k$. The case $k=1$ coincides
with Proposition \ref{p:1}. Suppose that the Proposition is valid for some
fixed $k$. Consider a set $\{f_1, \ldots, f_k, f_{k+1}\}$ of holomorphic
$V$-valued functions on $U$ having rank $k+1$. It is clear that the rank of
the set $\{f_1, \ldots, f_k\}$ equals $k$. Hence, one can find a set
$\{g_1, \ldots, g_k\}$ of holomorphic $V$-valued functions on $U$,
satisfying the requirement of the Proposition.

Taking into account Proposition \ref{p:1}, represent $f_{k+1}$ in the form
\[
f_{k+1} = g \,  d,
\]
where the holomorphic $V$-valued function $g$ has no zeros, and consider
the holomorphic $\bigwedge^{k+1}(V)$-valued function $g_1 \wedge \cdots
\wedge g_k \wedge g$. If this function has no zeros, then we take $g$ as
$g_{k+1}$ and get the set of functions $\{g_1, \ldots, g_k, g_{k+1}\}$
which satisfies the requirement of the Proposition. 

Suppose now that the function $g_1 \wedge \cdots \wedge g_k \wedge g$ has
zeros at the points of the set $\{c_\iota\}_{\iota \in {\mathcal I}}$ and
the zero at $c_\iota$ is of order $\mu_\iota$. Since $g_1(c_\iota) \wedge
\cdots \wedge g_k(c_\iota) \wedge g(c_\iota) = 0$ and the vectors
$g_1(c_\iota)$, $\ldots$, $g_k(c_\iota)$ are linearly independent one has
\[
g(c_\iota) = \sum_{\beta=1}^k g_\beta (c_\iota) \, b^{\beta}_\iota,
\]
where $b^{\beta}_{\iota}$ are some complex numbers. Using the
generalisation of the Weierstrass theorem and the generalisation of the
Mittag--Leffler theorem  \cite{Mar77}, one can show that there exist
functions $b^{\beta}$, $\beta = 1, \ldots, k$, which are holomorphic on $U$
and have the property
\[
b^{\beta}(c_\iota) = b^{\beta}_\iota.
\]
Define the holomorphic $V$-valued function
\[
h = g - \sum_{\beta=1}^k g_\beta \, b^{\beta}.
\]
It is clear that this function has zeros only at the points of the set
$\{c_\iota\}$. Moreover, since
\[
g_1 \wedge \cdots \wedge g_k \wedge \eta = g_1 \wedge \cdots \wedge g_k
\wedge g,
\]
the order of zero of $h$ at the point $c_\iota$ is less or equal to the
order of zero of $g_1 \wedge \cdots \wedge g_k \wedge g$. Represent $h$ in
the form
\[
h = g' \, d',
\]
where the holomorphic $V$-valued function $g'$ has no zeros. Then the
equality
\[
g_1 \wedge \cdots \wedge g_k \wedge g = g_1 \wedge \cdots \wedge g_k \wedge
g' \, d'
\]
implies that for any $\iota \in \mathcal I$ either the point $c_\iota$ is
not a zero of the function $g_1 \wedge \cdots \wedge g_k \wedge g'$, or it
is a zero of it but of order less than $\mu_\iota$. If the function $g_1
\wedge \cdots \wedge g_k \wedge g'$ has no zeros, we take $g'$ as $g_{k+1}$
and from the equality
\[
f_{k+1} = \left( \sum_{\beta=1}^k g_\beta \, b^{\beta} + g' \, d' \right) d
\]
conclude that the set $\{g_1, \ldots, g_k, g_{k+1}\}$ satisfies the
requirement of the Proposition. If it is again not the case, we repeat
along the lines of this paragraph using the function $g'$ instead of $g$.

Since the order of a zero is always finite, after a finite number of steps
we come to the set $\{g_1, \ldots, g_k, g_{k+1}\}$ satisfying the
requirement of the Proposition.
\end{proof}

From the proof of Proposition \ref{p:2} it follows that the functions
$g_1$, $\ldots$, $g_k$ can be chosen in such a way that the matrix $D =
(d^\beta_\alpha)$ is upper triangular. Moreover, if the set $\{f_1, \ldots,
f_l\}$, $l \le k$, is of constant rank $l$, then we can choose $g_1 = f_1$,
$\ldots$, $g_l = f_l$.

\begin{prop} \label{p:3}
Let $\{f_1, \ldots, f_k\}$, $k < n$, be a set of holomorphic $V$-valued
functions on $U$ of constant rank $k$. Then there exists a set $\{f_{k+1},
\ldots, f_n\}$ of holomorphic $V$-valued functions on $U$ such that the set
$\{f_1, \ldots, f_k, f_{k+1}, \ldots, f_n\}$ is of constant rank $n$ on
$U$.
\end{prop}

\begin{proof}
Let $c$ be some point of $U$. Denote $e_1 = f_1(c)$, $\ldots$, $e_k =
f_k(c)$. The vectors $e_1$, $\ldots$, $e_k$ are linearly independent, and
one can find vectors $e_{k+1}$, $\ldots$, $e_n$ such that $\{e_1, \ldots,
e_k\}$ is a basis of $V$. Consider the set $\{f_1, \ldots, f_k, e_{k+1}\}$
of holomorphic $V$-valued functions on $U$. It is clear that the rank of
this set equals $k+1$. Applying arguments used in the proof of Proposition
\ref{p:2}, one can show that there is a holomorphic $V$-valued function
$f_{k+1}$ on $U$ such that the set $\{f_1, \ldots, f_k,  f_{k+1}\}$ is of
constant rank $k+1$ on $U$. Using such a construction repeatedly we find a
set $\{f_1, \ldots, f_n\}$ satisfying the requirement of the Proposition.
\end{proof}

\begin{prop} \label{p:4}
Let $\{f_1, \ldots, f_n\}$ be a set of holomorphic $V$-valued functions on
$U$ of constant rank $n$. Then there exists a unique set $\{f^1, \ldots,
f^n\}$ of holomorphic $V^*$-valued functions on $U$ such that
\[
f^i(f_j) = \delta^i_j, \qquad i,j = 1, \ldots, n.
\]
\end{prop}

\begin{proof}
Let $\{e_i\}$ be a basis of $V$, and $\{e^i\}$ be the dual basis of $V^*$.
One has
\[
f_i = \sum_{j=1}^n e_j \, b^j_i,
\]
where $b^j_i$, $i, j = 1, \ldots, n$, are holomorphic functions on $U$.
Since for any $c \in U$ the vectors $f_1(c), \ldots, f_n(c)$ are linearly
independent, the matrix $(b^j_i(c))$ is nondegenerate. Hence, there are
holomorphic functions $d^i_j$, $i, j = 1, \ldots, n$, on $U$ such that for
any $c \in U$ one has
\[
\sum_{m = 1}^n d^i_m(c) \, b^m_j(c) = \delta^i_j, \qquad i, j = 1, \ldots,
n. 
\]
It is easy to get convinced that the holomorphic $V^*$-valued functions 
\[
f^i = \sum_{j=1}^n d^i_j \, e^j, \qquad i = 1, \ldots, n,
\]
satisfy the requirement of the Proposition.
\end{proof}

\begin{prop} \label{p:5}
Let $\{f_1, \ldots, f_k\}$ be a set of holomorphic $V$-valued functions on
$U$ of constant rank $k$, and $g$ be a holomorphic $V$-valued function on
$U$ such that $f_1 \wedge \cdots \wedge f_k \wedge g \equiv 0$. Then there
is a unique representation
\[
g = \sum _{\alpha=1}^k f_\alpha \, b^\alpha,
\]
where $b^\alpha$, $\alpha=1, \ldots, k$, are holomorphic functions on $U$.
\end{prop}

\begin{proof}
For any $c \in U$ the vectors $f_1(c)$, $\ldots$, $f_k(c)$ are linearly
independent, while the vectors $f_1(c)$, $\ldots$, $f_k(c)$, $g(c)$ are
linearly dependent, therefore one has
\[
g = \sum_{\alpha=1}^k f_\alpha \, b^\alpha,
\]
where $b^\alpha$, $\alpha = 1, \ldots, k$, are complex functions on $U$.
Due to Proposition \ref{p:3} one can construct holomorphic $V$-valued
functions $f_{k+1}$, $\ldots$, $f_n$ such that the set $\{f_1, \ldots,
f_n\}$ is of constant rank $n$. Let $f^i$, $i = 1, \ldots, n$, be the
holomorphic $V^*$-valued functions satisfying the requirement of
Proposition \ref{p:4}. Using these functions, one gets
\[
b^\alpha = f^\alpha (g), \qquad \alpha = 1, \ldots, k.
\]
Thus, the functions $b^\alpha$, $\alpha = 1, \ldots, k$, are holomorphic.
\end{proof}

\begin{prop} \label{t:1}
Let the set $\{f_1, \ldots, f_k\}$ of holomorphic $V$-valued functions on
$U$ be of rank $l \le k$. Then there is a set $\{g_1, \ldots, g_l\}$ of
holomorphic $V$-valued functions on $U$ of constant rank $l$ such that
there is valid the representation
\begin{equation}
f_\alpha = \sum_{\beta=1}^l g_\beta \, d^\beta_\alpha, \qquad \alpha = 1,
\ldots, k, \label{6}
\end{equation}
where $d^\beta_\alpha$, $\alpha = 1, \ldots, k$, $\beta = 1, \ldots, l$,
are holomorphic functions on $U$.

If $\{g'_1, \ldots, g'_l\}$ is another set of holomorphic $V$-valued
functions on $U$ of constant rank $l$ and for some holomorphic functions
$d'{}^{\beta}_\alpha$ one has
\[
f_\alpha = \sum_{\beta=1}^l g'_\beta \, d'{}^{\beta}_\alpha, \qquad \alpha
= 1, \ldots, k,
\]
then for some holomorphic functions $c^\beta_\alpha$, $\alpha, \beta = 1,
\ldots, l$, one has
\begin{equation}
g'_\alpha = \sum_{\beta=1}^l g_\beta \, c^\beta_\alpha, \qquad \alpha = 1,
\ldots, l. \label{5}
\end{equation}
\end{prop}

\begin{proof}
Without loss of generality we can assume that the rank of the set $\{f_1,
\ldots, f_l\}$ equals $l$. Thus, in accordance with Proposition \ref{p:2}
one can find a set $\{g_1, \ldots, g_l\}$ of holomorphic $V$-valued
functions on $U$ such that the set $\{g_1, \ldots, g_l\}$ is of constant
rank $l$ and (\ref{6}) is valid for $\alpha = 1, \ldots, l$. It is not
difficult to understand that for any $\alpha$ such that $l < \alpha \le k$
one has $g_1 \wedge \cdots \wedge g_l \wedge f_\alpha \equiv 0$. Thus, due
to Proposition \ref{p:5} representation (\ref{6}) is valid for all values
of $\alpha$.

To prove the second part of the Proposition, we suppose that the set
$\{f_1, \ldots, f_l\}$ is again of rank $l$. Then one has
\begin{eqnarray}
&&f_1 \wedge \cdots \wedge f_l = g_1 \wedge \cdots \wedge g_l \, \det
\widetilde D, \label{3} \\
&&f_1 \wedge \cdots \wedge f_l = g'_1 \wedge \cdots \wedge g'_l \, \det
\widetilde D', \label{4}
\end{eqnarray}
where $\widetilde D$ and $\widetilde D'$ are the holomorphic $\Mat(l,
\C)$-valued functions formed by the functions $d^\beta_\alpha$ and
$d'{}^\beta_\alpha$, $\alpha, \beta = 1, \ldots, l$, respectively. From
(\ref{4}) it follows that
\[
f_1 \wedge \cdots \wedge f_l \wedge g'_\alpha \equiv 0, \qquad \alpha = 1,
\ldots, l.
\]
Using this relation and taking into account (\ref{3}), one concludes that
\[
g_1 \wedge \cdots \wedge g_l \wedge g'_\alpha \equiv 0, \qquad \alpha = 1,
\ldots, l.
\]
Now Proposition \ref{p:5} implies that representation (\ref{5}) is valid.
\end{proof}

\begin{prop} \label{p:6}
Let $\{g_1, \ldots, g_k\}$ be a set of $V$-valued holomorphic functions on
$U$ of constant rank $k$, and $\{f_1, \ldots, f_l\}$ be a set of $V$-valued
holomorphic functions on $U$ such that the rank of the set $\{g_1, \ldots,
g_k, f_1, \ldots, f_l\}$ equals $k+m$, $m \le l$. Then there exist a set
$\{g_{k+1}, \ldots, g_{k+m}\}$ of $V$-valued holomorphic functions on $U$
such that the set $\{g_1, \ldots, g_k, g_{k+1}, \ldots, g_{k+m}\}$ is of
constant rank $k+m$ and there is valid the representation
\[
f_\alpha = \sum_{\beta=1}^{k+m} g_\beta d^\beta_\alpha, \qquad \alpha = 1,
\ldots, l,
\]
where $d^\beta_\alpha$, $\alpha =1, \ldots, l$, $\beta = 1, \ldots, k+m$,
are holomorphic functions on $U$.
\end{prop}

\begin{proof} 
The validity of the Proposition follows from the proof of Proposition
\ref{p:2}, from the discussion given after that proof, and from the proof
of Proposition \ref{t:1}.
\end{proof}

\subsection{Construction of Frenet frame}

Let $M$ be a Riemann surface, and $U$ be an open subset of $M$. The {\it
rank} of a holomorphic $\Mat(n, k; \C)$-valued function $\sm$ on $U$ is
defined as
\[
\rank \sm = \max_{p \in U} \rank \sm(p).
\]
We say that a holomorphic $\Mat (n, k; \C)$-valued function $\sm$ on $U$ is
of {\it constant rank} if
\[
\rank \sm(p)= \rank \sm
\]
for each $p \in U$. A $\mnk$-valued holomorphic function $\sm$ on $U$
generates the set $\{f_1, \ldots, f_k\}$ of $\C^n$-valued holomorphic
functions on $U$ determined by the columns of $\sm$. A $\mnk$-valued
holomorphic function $\sm$ on $U$ has rank $l$ if and only if the
corresponding set $\{f_1, \ldots, f_k\}$ is of rank $l$. A $\mnk$-valued
holomorphic function $\sm$ on $U$ is of constant rank if and only if the
corresponding set $\{f_1, \ldots, f_k\}$ is of constant rank. Having these
remarks in mind, we use the results of Section \ref{s:1} for $\mnk$-valued
holomorphic functions.

Consider again a holomorphic curve $\m : M \to \gkn$ and some holomorphic
lift $\sm: U \subset M \to \mnkt$ of $\m$. The function $\sm$ considered as
a $\mnk$-valued function is of constant rank $k$.

Since, in general, the function $\sm$ will be the first member of some
finite sequence of holomorphic matrix-valued functions, it is convenient to
denote $\sm_0 = \sm$ and $k_0 = k$.
Consider the holomorphic $\Mat(n, 2 k_0; \C)$-valued function $(\sm_0\;
\partial \sm_0/ \partial z)$ on $U$. This function is of rank $k_0 + k_1$,
where $0 \le k_1 \le k_0$. In the case $k_1 = 0$ due to Proposition
\ref{p:6} one has
\begin{equation}
\partial_- \sm_0 = \sm_0 \, B_{00}, \label{7}
\end{equation}
where $B_{00}$ is a holomorphic $\Mat(k_0; \C)$-valued function on $U$.
Therefore, all the derivatives of $\sm_0$ over $z$ can be expressed via
$\sm_0$ by relation similar to (\ref{7}). If $k_1 > 0$, then again due to
Proposition \ref{p:6} there exists a holomorphic $\Mat(n, k_1; \C)$-valued
function $\sm_1$ of constant rank $k_1$ such that
\begin{equation}
\partial_- \sm_0 = \sm_0 \, B_{00} + \sm_1 \, B_{10}, \label{10}
\end{equation}
where $B_{00}$ is a holomorphic $\Mat (k_0; \C)$-valued function and
$B_{10}$ is a holomorphic $\Mat(k_1, k_0; \C)$-valued function. The linear
subspaces of $\C^n$ spanned by the columns of the matrix $(\sm_0(p) \,
\sm_1(p))$, $p \in U$, do not depend on the choice of functions $\sm_0$ and
$\sm_1$. Actually, they are determined only by the mapping $\m$. Using
these subspaces, we define a holomorphic subbundle $\underline{\m}_1$ of
$\underline{\C}^n$ called the {\it first osculating space} of $\m$, which
generates a holomorphic curve  ${\m}_1 : M \to G^{k_0+k_1}(\C^n)$ called
the {\it first associated curve} of $\m$. 

Consider the holomorphic $\Mat(n, k_0+2k_1; \C)$-valued function $(\sm_0 \,
\sm_1 \, \partial \sm_1 / \partial z)$. It is of rank $k_0 + k_1 + k_2$,
where $0 \le k_2 \le k_1$. Again, if $k_2 = 0$, one concludes that the
derivatives of $\sm_0$ and $\sm_1$ over $z$ can be expressed via $\sm_0$
and $\sm_1$. If $k_2 \ne 0$ one finds a holomorphic $\Mat(n, k_2;
\C)$-valued function $\sm_2$ of constant rank $k_2$ such that
\[
\partial_- \sm_1 = \sm_0 \, B_{01} + \sm_1 \, B_{11} + \sm_2 \, B_{21}.
\]
After a finite number of steps we end up with the functions $\sm_0, \ldots,
\sm_t$ satisfying the relations
\begin{equation}
\partial_- \sm_a = \sum_{b=1}^{a+1} \sm_b \, B_{ba}, \label{2}
\end{equation}
where $B_{t+1, t} \equiv 0$. Thus, all the derivatives of $\sm_0$,
$\ldots$, $\sm_t$ over $z$ can be expressed via $\sm_0$, $\ldots$, $\sm_t$.
Each step of the construction gives us the corresponding subbundle
$\underline{\m}_a$, $a = 1, \ldots, t$, of $\underline \C^n$ and we get the
flag of the subbundles
\[
\underline{\m}_0 \subset \underline{\m}_1 \subset \cdots \subset
\underline{\m}_t.
\]

\begin{prop}
For any $p \in M$ the fiber $\underline \m_{tp}$ coincide with the linear
subspace $V$ defined in Section \ref{s:2}.
\end{prop}

\begin{proof}
The independence of $\underline \m_{tp}$ of $p$ can be proved along the
lines of the proof of Proposition \ref{t:2}. So we have just one subspace
of $\C^n$ which we denote by $W'$. It is quite clear that $W \subset W'$.
Suppose that $W$ is a proper subspace of $W'$. In this case there exists an
element $w \in \C^{n*}$ such that $w(v) = 0$ for each $v \in W$, but $w(u)
\ne 0$ for some element of $W'$. Identifying $\C^n$ with the space of
complex $n \times 1$ matrices and $\C^{n*}$ with the space of complex $1
\times n$ matrices, we can write
\[
w(v) = w v.
\]
The condition $w(v) = 0$ for each $v \in W$ is equivalent to the
requirement
\[
w \, \partial_-^l \sm = 0
\]
for $l = 0, 1, \ldots$. In particular, one has
\[
w \, \sm = w \, \sm_0 = 0. 
\]
Further, using (\ref{10}) one obtains
\[
w \, \partial_- \sm_0 = w \, \sm_1 \, B_{10} = 0.
\]
Since the $\Mat(k_1, k_0; \C)$-valued function $B_{10}$ is of rank $k_1$,
one has 
\[
w \, \sm_1 = 0.
\]
Similarly one shows that $w \, \sm_a = 0$ for all $a = 1, \ldots, t$, but
this contradicts to our supposition. Thus, $W' = W$.
\end{proof}

For a linearly full mapping $\m$ one has $W = \C^n$, therefore,
$\underline{\m}_t = \underline{\C}^n$. In any case the $t$th associated
curve is trivial.

Starting from this point, we assume that the linear space  $\C^n$ is
endowed with a positive definite hermitian scalar product 
\[
(v, u) = \sum_{i,j=1}^n \overline{v^i} \,h_{\bar \imath j} \, u^j,
\]
that in matrix notation can be written as
\[
(v, u) = v^\dagger h \, u.
\]
Here and below we denote by $A^\dagger$ the hermitian conjugate of the
matrix $A$.

Now we proceed to the construction of what we call the Frenet frame
associated with the lift $\sm$ of the mapping $\m$. Consider a $\Mat(n,
\C)$-valued function
\begin{equation}
\Pi_0 = I_n - \sm_0^{} (\sm_0^\dagger \, h \, \sm_0^{})^{-1} \sm_0^\dagger
\, h, \label{15}
\end{equation}
where $I_n$ is the $n \times n$ unit matrix.
For any $p \in U$ the matrix $\Pi_0(p)$ is a matrix of the operator of the
orthogonal projection to the orthogonal complement of $\underline\m_{0p}$
in $\C^n$, with respect to the canonical basis of $\C^n$. From (\ref{10})
one immediately gets
\begin{equation}
\Pi_0 \, \partial_- \sm_0 = \Pi_0 \, \sm_1 \, B_{10}. \label{12}
\end{equation} 
Denoting
\begin{equation}
\fm_0 = \sm_0, \qquad \fm_1 = \Pi_0 \, \sm_1, \label{13}
\end{equation}
and
\[
\beta_0 = \fm_0^\dagger \, h \, \fm_0^{},
\]
we can write (\ref{12}) in the form
\begin{equation}
\partial_- \fm_0 = \fm^{}_0 \, \beta_0^{-1} \partial_- \beta^{}_0 + \fm_1
\, B_{10}, \label{16}
\end{equation}
where we used the relation
\[
\partial_- \beta_0 = \fm_0^\dagger \, h \, \partial_- \fm^{}_0
\]
which follows from the equality
\[
\partial_+ \fm_0 = 0.
\]
The subspaces spanned by the linear combinations of the columns of the
matrix $\fm_1(p)$ do not depend on the choice of the lift $\sm = \sm_0$ and
the function $\sm_1$. These subspaces generate the subbundle
$\underline\varphi_1$ of $\underline \C^n$. Actually it is an orthogonal
complement of the subbundle $\underline \fm_0 = \underline \m_0$ in
$\underline \m_1$. Thus, we have the following representation
\[
\underline \m_1  = \underline \fm_0 \oplus \underline \fm_1,
\]
where the sum in the right hand side is orthogonal. Note that the
orthogonality of the fibers of $\underline \fm_0$ and $\underline \fm_1$ is
equivalent to the validity of the equality
\[
\fm^\dagger_0 \, h \, \fm_1^{} = 0.
\]

Find now the equations satisfied by the function $\fm_1$. Define the
$\Mat(n, \C)$-valued function $\Pi_1$ by
\[
\Pi_1 = I_n - \fm_1^{} \beta_1^{-1} \fm_1^\dagger h,
\]
where
\[
\beta_1 = \fm_1^\dagger \, h \, \fm_1^{}.
\]
For any $p \in U$ the matrix $\Pi_1(p)$ is the matrix of the operator of
the orthogonal projection to the orthogonal complement of $\underline
\fm_{1p}$ in $\C^n$. Since the subspaces $\underline \fm_{0p}$ and
$\underline \fm_{1p}$ are orthogonal, we have
\begin{equation}
\Pi_0 \, \Pi_1 = \Pi_1 \, \Pi_0. \label{14}
\end{equation} 
Using relations (\ref{13}), (\ref{14}) and the equality
\[
\partial_- \sm_1 = \sm_0 \, B_{01} + \sm_1 \, B_{11} + \sm_{2} \, B_{21},
\]
one gets
\begin{equation}
\Pi_2 \, \partial_- \fm_1 = \Pi_1 \, \partial_- \Pi_0 \, \sm_1 + \Pi_1 \,
\Pi_0 \, \sm_2 \, B_{21}. \label{17}
\end{equation}
It follows from (\ref{15}) and (\ref{16}) that
\[
\partial_- \Pi_0 = - \fm^{}_1 \, B^{}_{10} \, \beta_0^{-1} \fm_0^\dagger \,
h.
\]
Therefore, denoting
\[
\fm_2 = \Pi_1 \, \Pi_0 \, \sm_2,
\]
we rewrite (\ref{17}) as
\begin{equation}
\partial_- \fm_1 = \fm^{}_1 \, \beta_1^{-1} \fm_1^\dagger \, h \,
\partial_- \fm^{}_1 + \fm_2 \, B_{21}. \label{18}
\end{equation}
Similarly one gets
\[
\partial_+ \Pi_0 = h^{-1} \left( \partial_- \Pi_0 \right)^\dagger h =
\fm^{}_0 \, \beta_0^{-1} D^{}_{01} \fm_1^\dagger \, h,
\]
where
\[
D_{01} = - B_{10}^\dagger
\]
is an antiholomorphic $\Mat(k_1, k_2; \C)$-valued function. Hence,
\[
\partial_+ \fm_1 = \fm^{}_0 \, \beta_0^{-1} D^{}_{01} \, \beta_1.
\]
This equation, in particular, implies that
\[
\partial_- \beta_1 = \fm_1^\dagger \, h \, \partial_- \fm_1.
\]
Therefore, equation (\ref{18}) can be written as
\[
\partial_- \fm_1 = \fm^{}_1 \, \beta_1^{-1} \partial_- \beta^{}_1 + \fm_2
\, B_{21}. 
\]
As before, the columns of the matrices $\fm_2(p)$ generate a subbundle
$\underline \fm_2$ and we have the orthogonal decomposition
\[
\underline \m_2 = \underline \fm_0 \oplus \underline \fm_1 \oplus
\underline \fm_2.
\]

In general, defining inductively the functions
\[
\beta_a = \fm^\dagger_a \, h \, \fm^{}_a,
\]
the projectors
\[
\Pi_a = I_n - \fm_a^{} \beta_a^{-1} \fm_a^\dagger \, h
\]
and the functions
\[
\fm_{a+1} = \Pi_a \cdots \Pi_0 \, \sm_{a+1},
\]
we get the orthogonal decompositions
\[
\underline \m_a = \bigoplus_{b=0}^a \underline \fm_b.
\]
Note that the sequence of subbundles $\underline \fm_1, \ldots, \underline
\fm_t$ and the sequence of the mappings from $M$ to the corresponding
Grassmann manifolds are partial cases of the so called harmonic sequences
of subbundles and mappings intensively used to classify and construct
harmonic mappings from a Riemann surface to a Grassmann manifold. 

\begin{theor}
The functions $\fm_a$, $a = 0, \ldots, t$, satisfy the equations
\begin{eqnarray}
& \displaystyle \partial_- \fm_a = \fm^{}_a \, \beta_a^{-1} \partial_-
\beta^{}_a + \fm_{a+1} \, B_{a+1, a}, \label{19} \\
& \displaystyle \partial_+ \fm_a = \fm^{}_{a-1} \, \beta_{a-1}^{-1}
D^{}_{a-1,a} \, \beta^{}_a, \label{20}
\end{eqnarray}
where $B_{t+1,t} \equiv 0$, $D_{-1,0} \equiv 0$ and
\begin{equation}
D_{a-1,a} = - (B_{a, a-1})^\dagger. \label{2.1}
\end{equation}
\end{theor}

\begin{proof}
Assume that that we proved the validity of equations (\ref{19}), (\ref{20})
for all $a \le b < t$. Using these equations for $a=b$, one has
\begin{eqnarray*}
&& \partial_- \Pi_b = -\varphi_{b+1} \, B_{b+1, b} \, \beta_b^{-1}
\varphi_b^\dagger h + \varphi_b \, B_{b, b-1} \, \beta_{b-1}^{-1} \,
\varphi_{b-1}^\dagger h, \\
&& \partial_+ \Pi_b = -\varphi_{b-1} \, \beta_{b-1}^{-1} \, D_{b-1, b} \,
\varphi_b^\dagger h + \varphi_b \, \beta_b^{-1} D_{b, b+1} \,
\varphi_{b+1}^\dagger h. 
\end{eqnarray*}
Now differentiating the definition of $\fm_{b+1}$ over $z^-$ and $z^+$, we
can easily show that equations (\ref{19}), (\ref{20}) are valid for $a =
b+1$.
\end{proof}

Consider the $\Mat(n, \C)$-valued mapping $\varphi = (\varphi_0 \cdots
\varphi_t)$. For any $p \in U$ the matrix $\varphi(p)$ is nondegenerate.
Denote by $f_1, \ldots, f_n$ the $\C^n$-valued functions determined by the
columns of $\varphi$. For any $p \in U$ the vectors $f_1(p), \ldots,
f_n(p)$ are linearly independent and form a basis in $\C^n$. We call the
set of mappings $f_1, \ldots, f_n$, or the mapping $\varphi$ which
generates this set, the {\it Frenet frame} of $\psi$ associated with the
lift $\xi$.

\section{Connection with Toda systems}

\subsection{${\mathbb Z}$-graded Lie algebras}

To get the equations describing a Toda system \cite{LSa92, RSa94, RSa97a,
RSa97b} one starts with a complex Lie group $G$ whose Lie algebra
${\mathfrak g}$ is endowed with a ${\mathbb Z}$-gradation
\[
{\mathfrak g} = \bigoplus_{m \in {\mathbb Z}} {\mathfrak g}_m.
\]
Introduce the following subalgebras of ${\mathfrak g}$
\[
\tilde{\mathfrak h} = {\mathfrak g}_0, \qquad
\tilde{\mathfrak n}_- = \bigoplus_{m < 0} {\mathfrak g}_m, \qquad
\tilde{\mathfrak n}_+ = \bigoplus_{m > 0} {\mathfrak g}_m.
\]
and denote by $\tilde H$ and $\tilde N_\pm$ the connected Lie
subgroups of $G$ corresponding to the subalgebras $\tilde{\mathfrak
h}$ and $\tilde{\mathfrak n}_\pm$ respectively.

Suppose that $\widetilde H$ and $\widetilde N_\pm$ are closed subgroups of
$G$ and, moreover,
\begin{eqnarray*}
&\tilde H \cap \tilde N_\pm = \{e\}, \qquad \tilde N_-
\cap \tilde N_+ = \{e\}, \\
&\tilde N_- \cap \tilde H \tilde N_+ = \{e\}, \qquad
\tilde N_- \tilde H \cap \tilde N_+ = \{e\},
\end{eqnarray*}
where $e$ is the unit element of $G$. This is, in particular, true for the
finite-dimensional complex reductive Lie groups, see, for example,
\cite{Hum75}. Let the set $\tilde N_- \tilde H \tilde N_+$ be dense in $G$,
then for any element $g$ which belongs to $\tilde N_- \tilde H \tilde N_+$
one can write the
following unique decomposition
\begin{equation}
g = n_- h n_+^{-1}, \label{3.10}
\end{equation}
where $n_- \in \tilde N_-$, $h \in \tilde H$ and $n_+ \in
\tilde N_+$. This is again true for the finite-dimensional complex
reductive Lie groups.  Decomposition (\ref{3.10}) is called the {\it Gauss
decomposition}.

There is a simple classification of possible ${\mathbb Z}$-gradations
for complex semi\-simple Lie algebras, see, for example,
\cite{GOV94,RSa97b}. In this case for any ${\mathbb Z}$-gradation of a such
an algebra ${\mathfrak g}$, there exists a unique element $q \in {\mathfrak
g}$ which has the following property. An element $x \in {\mathfrak g}$
belongs to the subspace ${\mathfrak g}_m$ if and only if $[q, x] = m x$.
The element $q$ is called the {\it grading operator}.

Let ${\mathfrak g}$ be a semisimple Lie algebra of rank $r$. Denote by
$h_i$ and $i = 1, \ldots, r$, some set of the Cartan generators of
${\mathfrak g}$. For any set of $r$ nonnegative numbers $l_i$ the element
\begin{equation}
q = \sum_{i,j=1}^r h_i (k^{-1})_{ij} l_j, \label{2.5}
\end{equation}
where $k = (k_{ij})$ is the Cartan matrix of ${\mathfrak g}$, is the
grading operator of some ${\mathbb Z}$-gradation of ${\mathfrak g}$. The
numbers $l_i$ can be considered as  the labels assigned to the vertices of
the corresponding Dynkin diagram. The well-known canonical gradation of
${\mathfrak g}$ arises when one chooses all the numbers $l_i$ equal to 2.
Actually, if two sets of labels are connected by an automorphism of the
Dynkin diagram, we get two ${\mathbb Z}$-gradations connected by an
`external' automorphism of $\mathfrak g$. If it is not the case, then
different sets of labels give $\mathbb Z$-gradations which cannot be
connected by an automorphism of $\mathfrak g$.

If all the labels of the Dynkin diagram of a semisimple Lie algebra
${\mathfrak g}$ are different from zero, then the subgroup ${\mathfrak
g}_0$
coincides with the Cartan subalgebra ${\mathfrak h}$. In this case
the subgroup $\tilde H$ is abelian and we obtain the so called
{\it abelian Toda equations}. In all other cases the subgroup $\tilde
H$ is nonabelian and the corresponding Toda equations are called
{\it nonabelian}.

If we deal with a reductive Lie algebra, we choose as the grading operator
any
grading operator of its maximal semisimple subalgebra.

\subsection{Toda equations and their general solution}

Let $M$ be a simply connected complex one dimensional manifold. Consider a
reductive complex Lie group $G$ whose Lie algebra ${\mathfrak g}$ is
endowed with a ${\mathbb Z}$-gradation. Let $l$ be a positive integer, such
that the grading subspaces ${\mathfrak g}_m$ for $-l < m < 0$ and $0 <
m < l$ are trivial, and $c_-$ and $c_+$ be some fixed mappings from $M$ to
the subspaces ${\mathfrak g}_{-l}$ and ${\mathfrak g}_{+l}$, respectively,
such that
$c_-$ is holomorphic and $c_+$ is antiholomorphic. Restrict ourselves to
the case when $G$ is a matrix Lie group. In this case the Toda equations
are the matrix partial differential
equations of the form 
\[
\partial_+(\gamma^{-1} \partial_- \gamma) = [c_-, \gamma^{-1} c_+
\gamma],
\]
where $\gamma$ is a mapping from $M$ to $\tilde H$. These equations are the
zero curvature condition for the connection $\omega = \omega_- dz^- +
\omega_+ dz^+$ on the trivial principal $G$-bundle $M \times G$, where
\[
\omega_- = \gamma^{-1} \partial_- \gamma + c_-, \qquad \omega_+ =
\gamma^{-1} c_+ \gamma.
\]
Since $M$ is simply connected, then there exists a mapping $\fm: M \to G$
such that
\begin{equation}
\partial_- \fm = \fm (\gamma^{-1} \partial_- \gamma + c_-), \qquad
\partial_+ \fm = \fm \gamma^{-1} c_+ \gamma. \label{***}
\end{equation}

To obtain the general solution of Toda equations one can use the following
procedure \cite{LSa92,RSa94,RSa97a,RSa97b}. Choose some mappings
$\gamma_\pm$ from $M$ to $\tilde H$ such that $\gamma_-$ is holomorphic and
$\gamma_+$ is antiholomorphic.
Then integrate the equations
\begin{equation}
\mu^{-1}_\pm \partial_\pm \mu_\pm = \gamma_\pm c_\pm \gamma_\pm^{-1},
\qquad \partial_\mp \mu_\pm = 0. \label{2.3}
\end{equation}
The Gauss decomposition (\ref{3.10}) induces the corresponding
decomposition of mappings from $M$ to $G$. In particular, one obtains
\[
\mu^{-1}_+ \mu_- = \nu_- \eta \nu_+^{-1},
\]
where the mapping $\eta$ takes values in $\tilde H$, and the mappings
$\nu_\pm$ take values in $\tilde N_\pm$. It can be shown that the mapping
\[
\gamma = \gamma^{-1}_+ \eta \gamma_-
\]
satisfies the Toda equations, and any solution to these equations can be
obtained by the described procedure. Note that for the corresponding
mapping $\varphi$ one has the following expression \cite{RSa94, RSa97a,
RSa97b}
\begin{equation}
\varphi = g \mu_+ \nu_- \eta \gamma_- = g \mu_- \nu_+ \gamma_-, \label{3.1}
\end{equation}
where $g$ is an arbitrary constant element of $G$.  

We will be interested in the so called hermitian solutions of Toda
equations \cite{RSa94, RSa97b}\footnote{Actually in \cite{RSa94, RSa97b} we
use the term `real solutions'.}. Let the Lie algebra ${\mathfrak g}$ be
endowed with some antilinear antiautomorphism $\sigma$ which can be
extended to the antiholomorphic antiautomorphism $\Sigma$ of the Lie group
$G$. Suppose that
\[
\sigma(\mathfrak g_m) = \mathfrak g_{-m}, \qquad m \in \mathbb Z.
\]
In this case one has
\[
\Sigma(\widetilde H) = \widetilde H, \qquad  \Sigma(\widetilde N_\pm) =
\widetilde N_\mp.
\]
Assume that the mappings $c_\pm$ entering the Toda equations are subjected
to the condition
\[
\sigma(c_-) = -c_+.
\]
In this case, if $\gamma$ is a solution to Toda equations, then $\Sigma
\circ \gamma$ is also a solution of the same equations, and we can consider
solutions of Toda equations having the property
\begin{equation}
\Sigma \circ \gamma = \gamma. \label{23}
\end{equation}
To get such solutions using the general procedure of integration of the
Toda equations described above, we should start with the mappings
$\gamma_\pm$ which satisfy the relation
\[
\gamma_+ = \Sigma \circ \gamma_-^{-1}
\]
and choose the solutions of (\ref{2.3}) for which
\[
\mu_+ = \Sigma \circ \mu_-^{-1}.
\]
It can be shown that in such a way we get all the solutions of the Toda
equations satisfying (\ref{23}). We call such solutions {\it hermitian}.
The mapping $\varphi$ corresponding to a hermitian solution of Toda
equations satisfies the relation
\[
(\Sigma \circ \varphi) \Sigma(g) g \varphi = \gamma,
\]
where $g$ is the element of $G$ entering (\ref{3.1}).
 
\subsection{Nonabelian Toda equations associated with Lie group ${\rm
GL}(n, \C)$ and Frenet frames}

In this Section, generalising the consideration of \cite{RSa97c}, we give
the general form of nonabelian Toda equations based on the Lie group ${\rm
GL}(n, \C)$ and establish their connection to the equations of the Frenet
frame.

The Lie algebra $\mathfrak {gl}(n, {\mathbb C})$ of the Lie group ${\rm
GL}(n, \C)$ is reductive and can be represented as the direct product of
the simple Lie
algebra $\mathfrak {sl}(n, {\mathbb C})$ of rank $n-1$ and a one
dimensional Lie algebra isomorphic to $\mathfrak {gl}(1, {\mathbb C})$ and
composed by the $n
\times n$ complex matrices which are multiplies of the unit matrix.
Consider the general ${\mathbb Z}$-gradation of $\mathfrak {sl}(n, {\mathbb
C})$ arising when we choose the labels of the corresponding Dynkin diagram
as follows
\[
\begin{picture}(308,35) 
  \Line(2,4)-(12,4)
  \Line(38,4)-(98,4)
  \Line(124,4)-(184,4)
  \Line(210,4)-(270,4)
  \Line(296,4)-(306,4)
  \FilledCircle(2,4){2}[100]
  \FilledCircle(18,4){0.5}[0]
  \FilledCircle(25,4){0.5}[0]
  \FilledCircle(32,4){0.5}[0]  
  \FilledCircle(48,4){2}[100]
  \FilledCircle(68,4){2}[100]
  \FilledCircle(88,4){2}[100]
  \FilledCircle(104,4){0.5}[0]
  \FilledCircle(111,4){0.5}[0]
  \FilledCircle(118,4){0.5}[0]
  \FilledCircle(134,4){2}[100]
  \FilledCircle(154,4){2}[100]
  \FilledCircle(174,4){2}[100]
  \FilledCircle(190,4){0.5}[0]
  \FilledCircle(197,4){0.5}[0]
  \FilledCircle(204,4){0.5}[0]
  \FilledCircle(220,4){2}[100]
  \FilledCircle(240,4){2}[100]
  \FilledCircle(260,4){2}[100]
  \FilledCircle(276,4){0.5}[0]
  \FilledCircle(283,4){0.5}[0]
  \FilledCircle(290,4){0.5}[0]
  \FilledCircle(306,4){2}[100]
  \Text(2,9)[cb]{$\scriptstyle 0$}
  \Text(48,9)[cb]{$\scriptstyle 0$}
  \Text(68,9)[cb]{$\scriptstyle s_1$}
  \Text(88,9)[cb]{$\scriptstyle 0$}
  \Text(134,9)[cb]{$\scriptstyle 0$}
  \Text(154,9)[cb]{$\scriptstyle s_2$}
  \Text(174,9)[cb]{$\scriptstyle 0$}
  \Text(220,9)[cb]{$\scriptstyle 0$}
  \Text(240,9)[cb]{$\scriptstyle s_t$}
  \Text(260,9)[cb]{$\scriptstyle 0$}
  \Text(306,9)[cb]{$\scriptstyle 0$}
  \Text(25,14)[cb]{$\overbrace{\hbox to 50\unitlength{\hfil}}^{k_0-1}$}
  \Text(111,14)[cb]{$\overbrace{\hbox to 50\unitlength{\hfil}}^{k_1-1}$}
  \Text(283,14)[cb]{$\overbrace{\hbox to 50\unitlength{\hfil}}^{k_t-1}$}
\end{picture}
\]
It is convenient to take as a Cartan subalgebra for $\mathfrak {sl}(n,
{\mathbb C})$ the subalgebra consisting of diagonal $n \times n$
matrices with zero trace. Here the standard choice of the Cartan
generators is 
\[
(h_i)_{kl} = \delta_{ki} \delta_{li} - \delta_{k,i+1} \delta_{l,i+1},
\qquad i = 1, \ldots, n-1.
\]
With such a choice of Cartan generators, using (\ref{2.5}), we obtain the
following block matrix form of the grading operator
\[
q = \left( \begin{array}{ccccc}
\rho_0 I_{k_0} & 0 & \cdots & 0 & 0 \\
0 & \rho_1 I_{k_1} & \cdots & 0 & 0 \\
\vdots & \vdots & \ddots & \vdots & \vdots \\
0 & 0 & \cdots & \rho_{t-1} I_{k_{t-1}} & 0 \\
0 & 0 & \cdots & 0 & \rho_t I_{k_t}
\end{array} \right),
\]
where
\[
\rho_a = \frac{1}{n} \left( - \sum_{b=1}^{a} s_b \sum_{c=0}^{b-1} k_c +
\sum_{b=a+1}^t s_b \sum_{c=b}^t k_c \right).
\]
We use this grading operator to define a ${\mathbb Z}$-gradation of the Lie
algebra $\mathfrak {gl}(n, {\mathbb C})$. It is not difficult to describe
the arising grading subspaces of $\mathfrak{gl}(n, \C)$ and the relevant
subgroups of ${\rm GL}(n, \C)$. To this end we consider an element $x$ of
$\mathfrak{gl}(n, \C)$ as a $t \times t$ block matrix $(x_{ab})_{a,b=0}^t$,
where $x_{ab}$ is $k_a \times k_b$ matrix. For fixed $a \ne b$, the block
matrices $x$ having only the block $x_{ab}$ different from zero belong to
the grading subspace $\mathfrak g_m$ with
\[
m = \sum_{c=a+1}^b s_c, \quad a < b, \qquad m = - \sum_{c=b+1}^a s_c, \quad
a > b.
\]
The block-diagonal matrices form the grading subspace $\mathfrak g_0$.

Using the same block matrix representation for the elements of ${\rm GL}(n,
\C)$, we see that the subgroup $\widetilde H$ consists of all
block-diagonal nondegenerate matrices, and the subgroups $\widetilde N_-$
and $\widetilde N_+$ consist, respectively, of all block lower and upper
triangular matrices with unit matrices on the diagonal. 

The equations related to Frenet frame arise in the case when all integers
$s_a$, $a = 1, \ldots, t$, are equal to $2$. In this case the mappings
$c_\pm$ have to take values in the subspaces $\mathfrak g_{\pm 2}$ and
their general form is
\[
c_- = \left( \begin{array}{ccccc}
0 & 0 & \cdots & 0 & 0 \\
B_{10} & 0 & \cdots & 0 & 0 \\
\vdots & \vdots & \ddots & \vdots & \vdots \\
0 & 0 & \cdots & 0 & 0 \\
0 & 0 & \cdots & B_{t, t-1} & 0
\end{array} \right), \qquad
c_+ = \left( \begin{array}{ccccc}
0 & D_{01} & \cdots & 0 & 0 \\
0 & 0 & \cdots & 0 & 0 \\
\vdots & \vdots & \ddots & \vdots & \vdots \\
0 & 0 & \cdots & 0 & D_{t-1,t} \\
0 & 0 & \cdots & 0 & 0
\end{array} \right) 
\]
Parametrise the mapping $\gamma$ as
\[
\gamma = \left( \begin{array}{ccccc}
\beta_0 & 0 & \cdots & 0 & 0 \\
0 & \beta_1 & \cdots & 0 & 0 \\
\vdots & \vdots & \ddots & \vdots & \vdots \\
0 & 0 & \cdots & \beta_{t-1} & 0 \\
0 & 0 & \cdots & 0 & \beta_t
\end{array} \right).
\]
The corresponding Toda equations are
\begin{equation}
\partial_+(\beta_a^{-1} \partial_- \beta_a) = - \beta_a^{-1} D_{a, a+1}
\beta_{a+1} B_{a+1,a} + B_{a, a-1} \beta_{a-1} D_{a-1, a} \beta_a,
\label{3.2}
\end{equation}
where $B_{t+1, t} \equiv 0$ and $D_{-1,0} \equiv 0$.

Represent the mapping $\varphi$ entering (\ref{***}) in the block form
$\varphi = (\varphi_0 \cdots \varphi_t)$, where $\varphi_a$ is a $\Mat(n,
k_a; \C)$-valued function. Then it is clear that  relations (\ref{***})
literally coincide with equations (\ref{19}), (\ref{20}). Thus, the
integrability conditions of the equations satisfied by the Frenet frame are
the Toda equations (\ref{3.2}).

Not any solution of equations (\ref{3.2}) gives us a mapping $\varphi$
describing the Frenet frame. First of all we should consider the Toda
systems for which relations (\ref{2.1}) are satisfied. These relations are
equivalent to the equality $(c_-)^\dagger = - c_+$. Moreover, for the
Frenet frame one has $\beta_a^\dagger = \beta_a$, which is equivalent to
$\gamma^\dagger = \gamma$. Hence, we should consider only hermitian
solutions to the Toda equations corresponding to the antiholomorphic
antiautomorphism of ${\rm GL}(n, \C)$ generated by the ordinary hermitian
conjugation of matrices. Finally, for the Frenet frame one has
$\varphi^\dagger \, h \, \varphi = \gamma$. Therefore, constructing the
mapping $\varphi$ with the help of (\ref{3.1}) one should take the constant
element $g$ of $G$ satisfying the relation $g^\dagger g = h$.

\section{Geometrical interpretation of Frenet frames}

\subsection{Invariant metrics on principal fibre bundles}

Let $P \stackrel{\pi}{\to} N$ be a principal fibre $G$-bundle, and $g^*$ be
a metric on $P$ which is invariant with respect to the standard right
action of $G$ on $P$. The metric $g^*$ determines a connection on $P$ for
which the horizontal subspace $\mathcal H_p \subset T_p(P)$, $p \in P$, is
defined as the orthogonal complement to the vertical subspace $\mathcal V_p
\subset T_p(P)$ tangent to the fiber through $p$. For any vector field $X$
on $N$ there is a unique horizontal vector field $X^*$ on $P$ satisfying
the relation 
\begin{equation}
d \pi \circ X^* = X \circ \pi. \label{4.4}
\end{equation}
The vector field $X^*$ is called the {\it horizontal lift} of $X$. The
vector field $X^*$ is invariant with respect to the standard right action
of $G$ on $P$. Below we call a vector field on $P$ invariant with respect
to the standard right action of $G$ simply an invariant vector field. Any
invariant horizontal vector field on $P$ is the horizontal lift of a unique
vector field on $N$. By definition, for any vector field $X$ on $N$ and a
function $f$ on $N$ one has
\begin{equation}
\pi^* (X(f)) = X^* (\pi^*f). \label{4.2}
\end{equation}
Denote by $X^\perp$ the component of the vector field $X$ on $P$ orthogonal
to the fibers, that is actually the horizontal component of $X$. For any
vector fields $X$ and $Y$ on $N$ the relation
\begin{equation}
[X, Y]^* = [X^*, Y^*]^\perp \label{4.3}
\end{equation}
is valid.

The metric $g^*$ determines a metric $g$ on $N$ defined by the equality
\begin{equation}
\pi^*(g(X, Y)) = g^*(X^*, Y^*). \label{4.1}
\end{equation}

\begin{prop} \label{pr:4.1}
The Levi--Civita connections $\nabla^*$ and $\nabla$ associated with $g^*$
and $g$ satisfy the equality
\[
(\nabla_Y X)^* = (\nabla^*_{Y^*} X^*)^\perp.
\]
\end{prop}

\begin{proof}
The validity of the assertion of the Proposition follows immediately from
the Koszul formula
\begin{eqnarray*}
2 g(\nabla_Y X, Z) &=& Y(g(X,Z)) + X(g(Z,Y)) - Z(g(Y,X)) \\
&-& g(Y, [X, Z]) + g(X, [Z, Y]) + g(Z, [Y, X])
\end{eqnarray*}
after taking into account definition (\ref{4.1}) and relations (\ref{4.2}),
(\ref{4.3}).
\end{proof}

Let $\m$ be an immersion of a manifold $M$ to $N$. A mapping $X: N \to
T(M)$ is said to be a vector field on $\m$ if
\[
\pr_M \circ X = \psi, 
\]
where $\pr_M$ is the canonical projection of $T(M)$ to $M$. A vector field
$X'$ defined on some open subset of $M$ is called an extension of a vector
field $X$ on $\psi$ if
\[
X' \circ \psi = X.
\]
Here and below writing relations containing extensions of vector fields we
imply the appropriate restrictions of the domains. 

The covariant derivative $\nabla_Y X$ of the vector field $X$ on $\psi$
along the vector field $Y$ on $M$ is defined as
\begin{equation}
\nabla_Y X = \nabla_{(d\psi \circ Y)'} X' \circ \psi, \label{4.5}
\end{equation}
where $X'$ and $(d\psi \circ Y)'$ are some arbitrary extensions of the
vector fields $X$ and $d\psi \circ Y$.

Let now $\sm: U \subset M \to P$ be  a local lift of $\m$. The mapping
$\sm$ is an immersion of $U$ to $P$. We call a vector field $X$ on $\sm$
{\it vertical} if for any $p \in U$ the vector $X_p$ is tangent to the
fiber through $p$. Similarly, a vector field $X$ on $\sm$ is said to be
{\it horizontal} if for any $p \in U$ the tangent vector $X_p$ is
horizontal, or, in other words, orthogonal to the fibre through $\sm(p)$.

Denote the induced connection generated by $\m$ also by $\nabla$.

\begin{prop} \label{pr:4.3}
For any vector field $X$ on $\psi$ and a vector field $Y$ on $M$ the
relation
\[
\nabla_Y X = d \pi \circ \nabla^*_{(d \psi \circ Y)^{\prime *}} X^{\prime*}
\circ \xi
\]
is valid.
\end{prop}

\begin{proof}
From Proposition \ref{pr:4.1} it follows that
\[
(\nabla_{(d\psi \circ Y)'} X')^* = (\nabla^*_{(d\psi \circ Y)^{\prime *}}
X^{\prime *})^\perp. 
\]
Relation (\ref{4.4}) gives
\[
\nabla_{(d\psi \circ Y)'} X' \circ \pi = d \pi \circ \nabla^*_{(d\psi \circ
Y)^{\prime*}} X^{\prime *}. 
\]
The statement of the Proposition follows now from definition (\ref{4.5}).
\end{proof}

\begin{prop} \label{pr:4.2}
Let $X$ be a horizontal vector field on $\xi$. There is an extension $X'$
of $X$ such that $X'$ is horizontal and invariant.
\end{prop}

\begin{proof}
The vector field $d \pi \circ X$ is a vector field on $\psi$. Let $(d \pi
\circ X)'$ be an arbitrary extension of $d \pi \circ X$. Consider the
vector field
\begin{equation}
X' = (d \pi \circ X)^{\prime *}. \label{4.6}
\end{equation}
This vector field is, by definition, horizontal and invariant. Let us show
that it is an extension of $X$. Indeed, from (\ref{4.6}) one gets
\[
d \pi \circ X' \circ \xi = (d \pi \circ X)' \circ \psi = d \pi \circ X.
\]
Since, $X' \circ \xi$ and $X$ are horizontal vector fields, then $X' \circ
\xi = X$. Thus, $X'$ is an extension of $X$ which satisfies the
requirements of the Proposition. 
\end{proof}

Below by an extension of a vector field $X$ on $\xi$ we mean an extension
of $X$ which satisfies the requirements of Proposition \ref{pr:4.2}.

\begin{prop} \label{pr:4.4}
Let $X$ be a horizontal vector field on $\sm$, and $Y$ be a vector field on
$M$. For the vector field $Z = d\pi \circ X$ on $\psi$ one has
\[
\nabla_Y Z = d \pi \circ \nabla^*_{(d \xi \circ Y)^{\perp \prime}} X' \circ
\xi.
\]
\end{prop}

\begin{proof}
Since $X'$ is horizontal and invariant, there is a unique vector field $Z'$
on $N$ such that
\[
Z^{\prime *} = X'.
\]
It can be easily shown that $Z'$ is an extension of $Z$. From Proposition
\ref{pr:4.3} it follows that
\[
\nabla_Y Z = d \pi \circ \nabla^*_{(d \psi \circ Y)^{\prime *}} X' \circ
\xi.
\]
Further, the vector field $(d \xi \circ Y)^{\perp \prime}$ is horizontal
and invariant. Therefore, there is a unique vector field $Y'$ such that
\[
(d \xi \circ Y)^{\perp \prime} = Y^{\prime *}.
\]
It is not difficult to demonstrate that $Y'$ is an extension of the vector
field $d\psi \circ Y$. Thus, we see that the assertion of the Proposition
is true.
\end{proof}

\subsection{Construction of a local moving frame}

Let us return to the case of holomorphic curves in Grassmann manifolds.
Denote by $\scf^i_\alpha$, $i=1,\ldots,n$, $\alpha=1,\ldots,k$, the
standard coordinate functions on $\mnk$ and their restrictions to $\mnkt$.
A generic vector field on $\mnkt$ has a unique representation
\[
X = \sum_{i=1}^n \sum_{\alpha=1}^k \left( X^i_\alpha
\frac{\partial}{\partial \scf^i_\alpha} + X^{\bar \imath}_{\bar \alpha}
\frac{\partial}{\partial \bar \scf^{\bar \imath}_{\bar \alpha}} \right).
\]

Recall that $\mnkt$ is the principal ${\rm GL}(k, \C)$-bundle over $\gkn$.
The vertical subspaces are generated by the Killing vector fields
corresponding to the right action of the Lie group ${\rm GL}(k, \C)$ on
$\mnkt$. The vector fields 
\[
K_\alpha^\beta = \sum_{i=1}^n \scf^i_\alpha \frac{\partial}{\partial
\scf^i_\beta}, \qquad \alpha, \beta = 1, \ldots, k,
\]
and
\[
\bar K_{\bar \alpha}^{\bar \beta} = \sum_{i=1}^n \bar \scf^{\bar
\imath}_{\bar \alpha} \frac{\partial}{\partial \bar \scf^{\bar
\imath}_{\bar \beta}}, \qquad \alpha, \beta = 1, \ldots, k,
\]
form bases in the spaces of the holomorphic and antiholomorphic Killing
vector fields. A vector field $X$ on $\mnkt$ is invariant with respect to
the right action of ${\rm GL}(k, \C)$ if and only if
\[
[K_\alpha^\beta, X] = 0, \qquad [\mbar K_{\bar \alpha}^{\bar \beta}, X] =
0, \qquad \alpha, \beta = 1, \ldots, k,
\]
that is equivalent to
\begin{eqnarray*}
&\displaystyle \sum_{i=1}^n \scf^i_\alpha \frac{\partial
X^j_\gamma}{\partial \scf^i_\beta} = X^j_\alpha \delta^\beta_\gamma, \qquad
\sum_{i=1}^n \scf^i_\alpha \frac{\partial  X^{\bar \jmath}_{\bar
\gamma}}{\partial \scf^i_\beta} = 0, \\ 
&\displaystyle \sum_{i=1}^n \bar \scf^{\bar \imath}_{\bar \alpha}
\frac{\partial  X^j_\gamma}{\partial \bar \scf^{\bar \imath}_{\bar \beta}}
= 0, \qquad
\sum_{i=1}^n \bar \scf^{\bar \imath}_{\bar \alpha} \frac{\partial X^{\bar
\jmath}_{\bar \gamma}}{\partial \bar\scf^{\bar \imath}_{\bar \beta}} =
X^{\bar \jmath}_{\bar \alpha} \delta^{\bar \beta}_{\bar \gamma}.
\end{eqnarray*}

To construct an invariant metric on $\mnkt$ define first the $\Mat(k,
\C)$-valued function $\Delta = (\Delta_{\bar \alpha \beta})$ on $\mnkt$ by
\[
\Delta_{\bar \alpha \beta} =  \sum_{i,j=1}^n \mbar \scf^{\bar \imath}_{\bar
\alpha} h^{}_{\bar \imath j} \scf^j_\beta.
\]
It is quite clear that the metric $g^*$ on $\mnkt$ given by 
\[
g^* = \sum_{i,j=1}^n \sum_{\alpha, \beta=1}^k d\mbar \scf^{\bar
\imath}_{\bar \alpha} \otimes h^{}_{\bar \imath j} \, d \scf^j_\beta \,
(\Delta^{-1})^{\beta \bar \alpha} + \sum_{i,j=1}^n \sum_{\alpha, \beta=1}^k
h^{}_{\bar \jmath i} \, d\scf^i_\alpha \, (\Delta^{-1})^{\alpha \bar \beta}
\otimes d \mbar \scf^{\bar \imath}_{\bar \beta}
\]
is invariant with respect to the right action of ${\rm GL}(k, \C)$ on
$\mnkt$. Hence, we can define the metric $g$ on $\gkn$ generated by $g^*$.
This metric is proportional to the standard K\"ahler metric on $\gkn$.

\begin{prop} \label{pr:4.5}
For any horizontal vector fields $X$ and $Y$ on $\mnkt$ one has
\begin{eqnarray*}
(\nabla_Y X)^\perp &=& \sum_{i,j=1}^n \sum_{\alpha, \beta=1}^k Y^j_\beta
\left( \frac{\partial X^i_\alpha}{\partial \scf^j_\beta}
\frac{\partial}{\partial \scf^i_\alpha} + \frac{\partial X^{\bar
\imath}_{\bar \alpha}}{\partial \scf^j_\beta} \frac{\partial}{\partial \bar
\scf^{\bar \imath}_{\bar \alpha}} \right)^\perp \\
&+& \sum_{i,j=1}^n \sum_{\alpha, \beta=1}^k Y^{\bar \jmath}_{\bar \beta}
\left( \frac{\partial X^i_\alpha}{\partial \bar \scf^{\bar \jmath}_{\bar
\beta}} \frac{\partial}{\partial \scf^i_\alpha} + \frac{\partial X^{\bar
\imath}_{\bar \alpha}}{\partial \bar \scf^{\bar \jmath}_{\bar \beta}}
\frac{\partial}{\partial \bar \scf^{\bar \imath}_{\bar \alpha}}
\right)^\perp. 
\end{eqnarray*}
\end{prop}

\begin{proof}
A vector field $X$ on $\mnkt$ is orthogonal to the fibers if and only if
\[
\sum_{i,j=1}^n \mbar w^{\bar \jmath}_{\bar \beta} h_{\bar \jmath i}
X^i_\alpha = 0, \qquad \sum_{i,j=1}^n X^{\bar \imath}_{\bar \alpha} h_{\bar
\imath j} w^j_\beta = 0.
\]
The Christoffel symbols of the metric $g^*$ are
\begin{eqnarray*}
& \displaystyle \Gamma^{i \beta \gamma}_{\alpha j k} = - \frac{1}{2}
\sum_{m=1}^n \sum_{\delta=1}^k \left( \delta^i_j \delta^\gamma_\alpha
(\Delta^{-1})^{\beta \bar \delta} \mbar \scf^{\bar m}_{\bar \delta} h_{\bar
m k} + \delta^i_k \delta^\beta_\alpha (\Delta^{-1})^{\gamma \bar \delta}
\mbar \scf^{\bar m}_{\bar \delta} h_{\bar m j} \right), \\
& \displaystyle \Gamma^{i \bar \beta \gamma}_{\alpha \bar \jmath k} =
\Gamma^{i \gamma \bar \beta}_{\alpha k \bar \jmath} = - \frac{1}{2} \left(
\sum_{m=1}^n \delta^i_k (\Delta^{-1})^{\gamma \bar \beta} h_{\bar \jmath m}
\scf^m_\alpha - \sum_{\delta=1}^k \delta^\gamma_\alpha h_{\bar \jmath k}
\scf^i_\delta (\Delta^{-1})^{\delta \bar \beta} \right), \\
&\Gamma^{\bar \imath \bar \beta \bar \gamma}_{\bar \alpha \bar \jmath \bar
k} = \mbar{\Gamma^{i \beta \gamma}_{\alpha j k}}, \qquad \Gamma^{\bar
\imath \beta \bar \gamma}_{\bar \alpha j \bar k} = \Gamma^{\bar \imath \bar
\gamma \beta}_{\bar \alpha \bar k j} = \mbar{\Gamma^{i \bar \beta
\gamma}_{\alpha \bar \jmath k}}, \qquad \Gamma^{i \bar \beta \bar
\gamma}_{\alpha \bar \jmath \bar k} =  \mbar{\Gamma^{\bar \imath \beta
\gamma}_{\bar \alpha j k}} = 0. 
\end{eqnarray*}
Using the above relations one sees that the assertion of the Proposition is
true.
\end{proof}

For an arbitrary vector field $X$ on $\sm$ one has the representation
\[
X = \sum_{i=1}^n \sum_{\alpha=1}^k \left[ X^i_\alpha \left(
\frac{\partial}{\partial \scf^i_\alpha} \circ \xi \right) + X^{\bar
\imath}_{\bar \alpha} \left( \frac{\partial}{\partial \bar \scf^{\bar
\imath}_{\bar \alpha}} \circ \xi \right) \right] ,
\]
where $X^i_\alpha$ and $X^{\bar \imath}_{\bar \alpha}$ are some functions
on $M$.
 
Let $\varphi$ be the Frenet frame of $\psi$ associated with the lift $\xi$.
Define the set of vector fields $\widetilde E^\alpha_i$, $i = 1, \ldots,
n$, $\alpha = 1, \ldots, k$, on $\sm$ by 
\[
\widetilde E^\alpha_{i} = \sum_{j=1}^n \varphi^j_i \left(
\frac{\partial}{\partial \scf^j_\alpha} \circ \xi \right).
\]
The vector fields $\widetilde E^\alpha_\beta$, $\alpha, \beta = 1, \ldots,
k$,  are vertical, while the vector fields $\widetilde E^\alpha_\mu$,
$\alpha = 1, \ldots, k$, $\mu = k+1, \ldots n$, are horizontal.

The vector fields
\[
E^\alpha_\mu = d \pi \circ \widetilde E^\alpha_\mu, \qquad \alpha = 1,
\ldots, k, \quad \mu = k+1, \ldots, n,
\]
form a local basis in the space of vector fields on $\psi$. Therefore, one
has
\[
\nabla_{\partial_\pm} E^\alpha_\mu = \sum_{\nu=k+1}^n \sum_{\beta=1}^k
E^\beta_\nu (\Lambda_\pm)^{\nu \alpha}_{\beta \mu}.
\]
Rewrite the equations (\ref{19}) and (\ref{20}) as
\[
\partial_\pm f_i = \sum_{j=1}^n f_j (\lambda_\pm)^j_i
\]
where the $\C^n$-valued functions $f_1, \ldots, f_n$ are generated by the
columns of the mapping $\varphi$.

\begin{theor}
The connection coefficients $(\Lambda_\pm)^{\nu \alpha}_{\beta \mu}$ are
connected with the functions $(\lambda_\pm)^{\nu \mu}$ by the relations
\[
(\Lambda_-)^{\nu \alpha}_{\beta \mu} = \delta^\alpha_\beta
(\lambda_-)^\nu_\mu - (\beta_0^{-1} \partial_- \beta_0)^\alpha_\beta
\delta^\nu_\mu, \qquad (\Lambda_+)^{\nu \alpha}_{\beta \mu} =
\delta^\alpha_\beta (\lambda_+)^\nu_\mu.
\]
\end{theor}

\begin{proof}
The statement of the Theorem follows from Propositions \ref{pr:4.4} and
\ref{pr:4.5}.
\end{proof}

In conclusion of this Section we give the expressions for the metrics
induced on $M$ by the mappings $\psi_a$, $a = 0, \ldots, t-1$. From the
equations satisfied by the Frenet frame it follows immediately that for the
metric $g_a$ induced on $a$th associated curve one has
\[
g_a(\partial_+, \partial_-) = - \tr (\beta_a^{-1}D_{a, a+1} \beta_{a+1}
B_{a+1, a}) = \tr (\beta_a^{-1} B^\dagger_{a+1, a} \beta_{a+1} B_{a+1, a}). 
\]
Toda equations (\ref{3.2}) give
\[
\partial_+ \partial_- \ln \det \beta_a = g_a (\partial_+ , \partial_-) -
g_{a-1} (\partial_+, \partial_-).
\]
This equality implies that the function $\ln (\det \beta_0 \cdots \det
\beta_a)$ is a K\"ahler potential of the metric $g_a$. 

\section*{Acknowledgments}

It is a pleasure to thank Yu. I. Manin and M. V. Saveliev for their
interest to the work and for fruitful discussions. The author is grateful
to the Max--Planck--Institut f\"ur Mathematik in Bonn where this work was
completed for the warm hospitality and financial support. This work was
also partially supported by the Russian Foundation for Basic Research under
grant no.~98--01--00015, and by the INTAS grant no.~96--690.

\end{document}